\documentclass[11pt,twoside,reqno]{amsart}
\usepackage[utf8]{inputenc}
\usepackage[T1]{fontenc}
\usepackage{import}
\usepackage{newclude}
\usepackage{amsmath,amssymb,amsthm,mathtools}

\usepackage{anysize}
\usepackage{amscd}
\usepackage{stmaryrd}
\usepackage{wasysym}
\usepackage{mathrsfs}
\usepackage[scr=boondox]{mathalfa}
\usepackage{enumitem}
\usepackage{accents}
\usepackage{dsfont}
\usepackage{soul}
\usepackage{color}
\usepackage[all]{xy}
\usepackage{float}
\usepackage{bm}
\usepackage{makecell}
\usepackage{thmtools}      
\usepackage{setspace}
\usepackage{comment}
\usepackage{tikz}
\usepackage{tikz-cd}
\usepackage{mathdots}
\usepackage{graphicx}
\usepackage[myheadings]{fullpage}

\usepackage[hypertexnames=false,colorlinks,linkcolor=black,citecolor=black]{hyperref}
\usepackage[capitalize,nameinlink]{cleveref}

\numberwithin{equation}{section}
\newcommand\mtop{.95in}
\newcommand\mbottom{.95in}
\newcommand\mleft{1in}
\newcommand\mright{1in}
\usepackage[top = \mtop, bottom = \mbottom, left = \mleft, right=\mright]{geometry}

\DeclareMathOperator{\Mat}{Mat}

\newtheorem{thm}{Theorem}[section]

\newtheorem{prop}[thm]{Proposition}

\newtheorem{lemma}[thm]{Lemma}

\newtheorem{cor}[thm]{Corollary}
\theoremstyle{definition}
\newtheorem{defi}[thm]{Definition}

\newtheorem{rmk}[thm]{Remark}

\newcommand\reallywidehat[1]{%
\savestack{\tmpbox}{\stretchto{%
  \scaleto{%
    \scalerel*[\widthof{\ensuremath{#1}}]{\kern-.6pt\bigwedge\kern-.6pt}%
    {\rule[-\textheight/2]{1ex}{\textheight}}
  }{\textheight}%
}{0.5ex}}%
\stackon[1pt]{#1}{\tmpbox}%
}
\DeclareSymbolFont{bbold}{U}{bbold}{m}{n}
\DeclareSymbolFontAlphabet{\mathbbold}{bbold}

\makeatletter
\def\@tocline#1#2#3#4#5#6#7{\relax
  \ifnum #1>\c@tocdepth 
  \else
    \par \addpenalty\@secpenalty\addvspace{#2}%
    \begingroup \hyphenpenalty\@M
    \@ifempty{#4}{%
      \@tempdima\csname r@tocindent\number#1\endcsname\relax
    }{%
      \@tempdima#4\relax
    }%
    \parindent\z@ \leftskip#3\relax \advance\leftskip\@tempdima\relax
    \rightskip\@pnumwidth plus4em \parfillskip-\@pnumwidth
    #5\leavevmode\hskip-\@tempdima
      \ifcase #1
       \or\or \hskip 1em \or \hskip 2em \else \hskip 3em \fi%
      #6\nobreak\relax
    \hfill\hbox to\@pnumwidth{\@tocpagenum{#7}}\par
    \nobreak
    \endgroup
  \fi}
\makeatother

\newcommand{\Z}{\mathbb{Z}}

\newcommand{\F}{\mathbb{F}}

\newcommand{\E}{\mathbb{E}}

\renewcommand{\l}{\lambda}

\newcommand{\Y}{\mathbb{Y}}

\DeclareMathOperator{\Hom}{Hom}

\DeclareMathOperator{\Cok}{Cok}

\DeclareMathOperator{\Ext}{Ext}

\DeclareMathOperator{\Aut}{Aut}

\DeclareMathOperator{\GL}{GL}

\DeclareMathOperator{\Sur}{Sur}

\title{Universality of rational canonical form for random matrices over a finite field}
\author{Jiahe Shen}

\date{\today}

\begin{document}

\thanks{I thank Roger Van Peski for many helpful discussions and suggestions, and Jason Fulman for further
comments on the draft. I also thank the anonymous referee for pointing out that the arguments in a
previous version follow from existing results in Cheong-Yu \cite{cheong2023distribution}. This work is supported by NSF grant
DMS-2246576 and Simons Investigator grant 929852.
}

\maketitle

\begin{abstract}
In this note, we study the distribution of the rational canonical form of a random matrix over the finite field
$\F_p$, whose entries are independent and $\epsilon$-balanced with $\epsilon\in(0,1-1/p]$.
We show that, as the matrix size tends to infinity, the statistics converge to independent
Cohen-Lenstra distributions, demonstrating the universality of this asymptotic behavior. In particular, we recover,
as a special case, the uniform setting proved in the thesis of Fulman
\cite{fulman1997probability}.

Our proof uses the fact that the rational canonical form data of $A_n$ and the $\F_p[t]$-module structure
of the function field cokernel $\Cok(tI_n-A_n)$ determine each other uniquely. Consequently, our
question can be reformulated, equivalently, as the asymptotic distribution problem for this
cokernel, which has been established by Cheong-Yu \cite{cheong2023distribution}.
\end{abstract}



\tableofcontents

\section{Introduction}\label{sec: Introduction}

\subsection{Main results}

Let $p$ be a prime number, and $\F_p$ be the finite field of order $p$. Denote by $\Mat_n(\F_p)$ the set of $n\times n$ matrices with entries in $\F_p$, and $\GL_n(\F_p)\subset\Mat_n(\F_p)$ the subset of invertible matrices. This paper studies the universal behavior of a random matrix $A_n\in\Mat_{n}(\F_p)$ as $n$ goes to infinity.

Before presenting our results, let us first explain the simplest form of a finite field matrix, which is always the main focus of random matrix theory. Denote by 
$$\Y=\{\l=(\l_1,\l_2,\ldots):\l_1\ge\l_2\ge\ldots\ge0,\l_i\in\Z,\l_i=0 \text{ for all but finitely many }i\}$$ 
the set of integer partitions, and fix a matrix $A_n\in\Mat_n(\F_p)$. Let $S\subset\F_p[t]$ denote the subset of monic, irreducible, nonconstant polynomials. Assume that the irreducible decomposition of the characteristic polynomial of $A_n$ has the form 
\begin{equation}\label{eq: irreducible decomposition of characteristic polynomial}
\det(tI_n-A_n)=\prod_{i=1}^kf_i(t)^{e_i},\quad f_1,\ldots,f_k\in S,
\end{equation}
where $\sum_{i=1}^ke_i\cdot\deg f_i=n$. Then there exists $P\in\GL_n(\F_p)$ such that
\begin{equation}\label{eq: rational canonical form}
P^{-1}A_nP=\begin{pmatrix}
R_1 & 0 & 0 & \cdots & 0 \\
0 & R_2 & 0 & \cdots & 0 \\
\cdots & \cdots & \cdots & \cdots & \cdots \\
0 & 0 & 0 & \cdots & R_k \\
\end{pmatrix}.
\end{equation}
Here, each matrix $R_i\in\Mat_{e_i\cdot\deg f_i}(\F_p),1\le i\le k$, has the form
$$R_i=\begin{pmatrix}
C(f_i^{\lambda_1^{(f_i)}}) & 0 & 0 & \cdots & 0 \\
0 & C(f_i^{\lambda_2^{(f_i)}}) & 0 & \cdots & 0 \\
\cdots & \cdots & \cdots & \cdots & \cdots \\
0 & 0 & 0 & \cdots & C(f_i^{\lambda_{j_i}^{(f_i)}}) \\
\end{pmatrix},$$
where $\lambda^{(f_i)}=(\lambda_1^{(f_i)},\ldots,\lambda_{j_i}^{(f_i)})\in\Y$ is a partition with $\lambda_1^{(f_i)}+\cdots+\lambda_{j_i}^{(f_i)}=e_i$, and for all monic $f=t^r+a_{r-1}t^{r-1}+\cdots+a_0\in\F_p[t]$,  
$$C(f):=
\begin{pmatrix}
0 & 1 & 0 & \cdots & 0 \\
0 & 0 & 1 & \cdots & 0 \\
\vdots & \vdots & \ddots & \ddots & \vdots \\
0 & 0 & 0 & \ddots & 1 \\
-a_{0} & -a_{1} & -a_{2} & \cdots & -a_{r-1}
\end{pmatrix}$$
is the \emph{companion matrix}. Furthermore, $A_n$ is conjugate (under the action of $\GL_n(\F_p)$) to a unique matrix of the form in \eqref{eq: rational canonical form}, which we call the \emph{rational canonical form} of $A_n$. In this way, when $A_n\in\Mat_n(\F_p)$ is random, we associate to every $f\in S$ the random partition $\lambda^{(f)}=\lambda^{(f)}(A_n)$ (when $f$ does not appear in the irreducible decomposition, we let $\l^{(f)}$ be the zero partition).

Now we turn to our main results. From now on, $\epsilon\in(0,1-1/p]$ will be a fixed real number, and we say that a random variable $\xi\in\F_p$ is \emph{$\epsilon$-balanced} if $\mathbf{P}(\xi=a)\le 1-\epsilon$ for all $a\in\F_p$. The following theorem shows that, when we enlarge the entry distribution from the uniform case proved by Fulman 
\cite[Section 3]{fulman1997probability} (also see Fulman \cite[Section 2.1]{fulman2002random}) to the broader class of independent $\epsilon$-balanced distributions, the asymptotic behavior of the rational canonical forms remains unchanged.


\begin{thm}\label{thm: joint CL distribution}
For each $n$, let $A_n\in\Mat_n(\F_p)$ be a random matrix with independent entries that are $\epsilon$-balanced. For any tuple of distinct polynomials $f_1,\ldots,f_r\in S$, the distribution of the $r$-tuple of random partitions $(\l^{(f_1)},\ldots,\l^{(f_r)})$ converges pointwise to the product measure $\prod_{i=1}^r\mu_{f_i}$ on $\Y^r$ as $n$ goes to infinity, where for all $1\le i\le r$,
$$\mu_{f_i}(\lambda)=\frac{1}{\#\Aut_{f_i}(\lambda)}\prod_{j=1}^\infty(1-p^{-j\deg f_i}),\quad\forall \lambda=(\lambda_1,\lambda_2,\ldots)\in\Y.$$
Here $\Aut_{f_i}(\lambda)$ is the group of automorphisms of the $\F_p[t]$-module $\bigoplus_{j\ge 1}\F_p[t]/(f_i^{\lambda_j})$.
\end{thm}

To illustrate the content of \Cref{thm: joint CL distribution}, let us first consider the most classical case $r = 1$ with $f_1(t)=t$.
In this situation, the partition $\l^{(t)}$ records the sizes of the nilpotent Jordan blocks of $A_n$, and \Cref{thm: joint CL distribution} asserts that 
$\lambda^{(t)}$ converges in distribution to the measure $\mu_{t}$ on integer partitions. This measure was introduced by Cohen and
Lenstra \cite{cohen2006heuristics} in their study of the distribution of class groups of imaginary quadratic fields. In its full form, \Cref{thm: joint CL distribution} extends to general $r \ge 1$, where the random
partitions $(\lambda^{(f_1)},\ldots,\lambda^{(f_r)})$ become asymptotically independent,
each following its own Cohen-Lenstra distribution. One can refer to
\Cref{prop: decomposition and explicit expression of Aut} for the explicit expression and properties of $\mathrm{Aut}_{f_i}(\lambda)$.

As a direct consequence of \Cref{thm: joint CL distribution}, by considering natural statistics such as the size $|\lambda^{(f_i)}|$ of each random
partition, we obtain the following corollary, which makes the
asymptotic behavior more concrete.

\begin{cor}\label{cor: joint multiplicity}
Let $A_n,f_1,\ldots,f_r$ be the same as in \Cref{thm: joint CL distribution}. For all $1\le i\le r$, denote by $e_i=|\lambda^{(f_i)}|$ the power of $f_i$ in the characteristic polynomial of $A_n$. Then the distribution of the $r$-tuple of non-negative integers $(e_1,\ldots,e_r)$ converges pointwise to the product measure $\prod_{i=1}^r\nu_{f_i}$ on $\Z_{\ge 0}^r$ as $n$ goes to infinity, where for all $1\le i\le r$,
$$\nu_{f_i}(j)=p^{-j\deg f_i}\prod_{k\ge j+1}^\infty(1-p^{-k\deg f_i}),\quad\forall j\in\Z_{\ge 0}.$$
\end{cor}

Taken together, \Cref{thm: joint CL distribution} and \Cref{cor: joint multiplicity} situate our work within a broader line of research on asymptotic behavior in discrete random matrix models. For the $p$-adic or integer case, Cohen-Lenstra distribution often arises as the limiting distribution of cokernels of random matrices. Maples \cite[Theorem 1.1]{maples2013cokernels} proved that the Cohen-Lenstra distribution is universal for $\Cok(A_n)$, $A_n\in\Mat_n(\Z_p)$ when $n$ goes to infinity. This extends the Haar measure case previously established by Friedman-Washington \cite{friedman1989distribution}. In his setting, the entries are required to be i.i.d.\ and $\epsilon$-balanced. The identically distributed restriction was then removed by Wood \cite[Theorem 1.2]{wood2019random}, who introduced the moment method that transfers the problem to studying the number of surjections. Following this work, Nguyen-Wood \cite{nguyen2022random} extended the moment method framework to study universal statistics for integer matrices, which also agree with the distributions defined by Cohen and Lenstra. Nguyen-Van Peski \cite{nguyen2024universality} and Huang-Nguyen-Van Peski \cite{huang2025cohen} generalized the result to the joint distribution of matrix products, showing that universality persists under matrix multiplication. Related developments in the study of cokernel statistics over more general rings and structured models include work of Cheong-Yu \cite{cheong2023distribution}, Van Werde \cite{van2025cokernel}, and Cheong-Huang \cite{cheong2025cokernel}.

The finite field case has also seen substantial progress. Fulman developed a generating function method in \cite{fulman1999probabilistic,fulman2002random}, later extended in collaboration with Neumann \cite{fulman2005generating}, to study random matrices over finite fields, for instance the probability that a random matrix is separable, cyclic, semisimple, or regular. This approach has also been connected to Cohen-Lenstra distributions, see the discussion in Fulman \cite{fulman2013cohen}. More recently, Luh-Meehan-Nguyen \cite{luh2021some} applied combinatorial and Fourier analytical techniques to obtain further results including the rank distribution and the uniformity of normal vectors for matrices with i.i.d.\ and $\epsilon$-balanced entries. We also refer to Fulman-Goldstein \cite{fulman2015stein}, Bl\H{o}mer-Karp-Welzl \cite{blomer1997rank}, Cooper \cite{cooper2000distribution}, Kahn-Komlós \cite{kahn2001singularity}, Maples \cite{maples2010singularity}, Ferber-Jain-Sah-Sawhney \cite{ferber2023random} for related results on singularity in various finite field random matrix models.

Although we work with a finite field model, we do not follow the methods from the finite field
literature listed above. Instead, the crucial input is
\Cref{thm: rational canonical and smith normal form}: it asserts that the rational canonical form data of
$A_n$ uniquely determine, and are uniquely determined by, the isomorphism class of the
$\F_p[t]$-module
$$
\Cok(tI_n-A_n):=\F_p[t]^n/(tI_n-A_n)\F_p[t]^n.
$$
The identification above transfers our problem to the asymptotics of
$\Cok(tI_n-A_n)$, viewed as a random finite $\F_p[t]$-module. Then, based on the universal distribution derived by
Cheong-Yu \cite{cheong2023distribution}, our main results follow immediately.

The above approach does not rely on Fulman’s generating function method. It therefore provides a
new proof in the uniform setting, and it also applies to the more general class of
$\epsilon$-balanced random matrices. As far as we know, this is the first demonstration that the
limiting distribution of rational canonical forms is insensitive to the specific distribution of
the entries.

\begin{rmk}
Throughout this paper, we require the finite field to have a prime order. In order to generalize these results to $\F_{q}$ where $q>p$ is a power of $p$, we need stronger condition for the distribution of the entries that they do not concentrate in proper affine subfields. Once these obstacles for Fourier analysis (see Wood \cite[Lemma 2.2]{wood2019random} for the integer case) are avoided, one can prove similar results as in \Cref{thm: joint CL distribution} and \Cref{cor: joint multiplicity} following the same steps.
\end{rmk}

\subsection{Notations}

We use $[n]$ to denote $\{1,2,\ldots,n\}$. We write $\F_p[t]$ as the polynomial ring with coefficients in $\F_p$. Denote by $\#$ the order of a finite set. We write $\mathbf{P}$ for probability and $\E$ for expectations. Given a ring $R$, we write $\Sur_{R}(\cdot,\cdot)$ for the set of surjections of $R$-modules, and  $\Hom_R(\cdot,\cdot)$ for the set of homomorphisms of $R$-modules.

\subsection{Outline of the paper}

In \Cref{sec: Preliminaries}, we introduce some necessary background for this paper, mainly around modules and matrices over the ring $\F_p[t]$. In \Cref{sec: Proof of the main results}, we prove the results stated in \Cref{sec: Introduction}.

\section{Preliminaries}\label{sec: Preliminaries}


\begin{defi}
We denote by $\Y$ the set of \emph{partitions} $\lambda=(\lambda_1,\lambda_2,\ldots)$, which are (finite or infinite) sequences of nonnegative integers $\lambda_1\ge\lambda_2\ge\cdots$ that are eventually zero. We do not distinguish between two such sequences that differ only by a string of zeros at the end. The integers $\l_i>0$ are called the \emph{parts} of $\l$. Let $|\l| := \sum_{i\ge 1}\l_i,n(\lambda):=\sum_{i\ge 1} (i-1)\lambda_i$, and $m_k(\l):= \#\{i\mid \l_i = k\}$. 
\end{defi}




Now we can classify the isomorphism classes of finite $\F_p[t]$-modules.

\begin{prop}\label{prop: classification of finite module}
Every finite $\F_p[t]$-module $G$ is isomorphic to one of the form
$$G\cong\bigoplus_{i=1}^r\bigoplus_{j\ge 1}\left(\F_p[t]/(f_i^{\l^{(f_i)}_j})\right),\quad \l^{(f_1)}=(\l^{(f_1)}_1,\ldots,\l^{(f_1)}_{j_1}),\ldots,\l^{(f_k)}=(\l^{(f_r)}_1,\ldots,\l^{(f_k)}_{j_r})\in\Y.$$
Here $f_1,\ldots,f_r\in S$ are distinct polynomials, and for all $1\le i\le r$, $\bigoplus_{j\ge 1}\left(\F_p[t]/(f_i^{\l^{(f_i)}_j})\right)$ is the $f_r$-submodule of $G$. Furthermore, let 
$$P:=f_1^{k_1}\cdots f_r^{k_r}\in\F_p[t]$$ 
be monic, where $f_1,\ldots,f_r\in S$, and let $R:=\F_p[t]/(P)$ be the finite quotient ring. Then, every finite $R$-module $M$ is isomorphic to one of the form
$$M\cong\bigoplus_{i=1}^r\bigoplus_{j\ge 1}\left(\F_p[t]/(f_i^{\l^{(f_i)}_j})\right),\quad \l^{(f_1)}=(\l^{(f_1)}_1,\ldots,\l^{(f_1)}_{j_1}),\ldots,\l^{(f_k)}=(\l^{(f_r)}_1,\ldots,\l^{(f_r)}_{j_r})\in\Y,$$
such that $\l_1^{(f_1)}\le k_1,\ldots,\l_1^{(f_r)}\le k_r$.
\end{prop}

In general, for $f\in S$, we say the finite $\F_p[t]$-module $G$ has \emph{type} $\l$ at $f$ if its $f$-submodule is isomorphic to $\left(\F_p[t]/(f^{\l_1})\right)\bigoplus\cdots\bigoplus\left(\F_p[t]/(f^{\l_j})\right)$. The following proposition studies the group of automorphisms of $\F_p[t]$-modules and counts their cardinality.

\begin{prop}\label{prop: decomposition and explicit expression of Aut}
We have
\begin{equation}\label{eq: decomposition of automorphism}
\Aut\left(\bigoplus_{i=1}^k\bigoplus_{j\ge 1}(\F_p[t]/(f_i^{\l^{(f_i)}_j}))\right)=\prod_{i=1}^k\Aut_{f_i}(\lambda^{(f_i)}),
\end{equation}
where for all $1\le i\le k$, $\Aut_{f_i}(\lambda^{(f_i)})$ is the group of automorphisms of $\bigoplus_{j\ge 1}(\F_p[t]/(f_i^{\l^{(f_i)}_j}))$, and
\begin{equation}\label{eq: explicit expression of automorphism}
\#\Aut_{f_i}(\lambda^{(f_i)})=p^{\deg f_i\cdot(|\l^{(f_i)}|+2n(\l^{(f_i)}))}\prod_{j\ge 1}(p^{-\deg f_i};p^{-\deg f_i})_{m_j(\l^{(f_i)})}.
\end{equation}
Here, the notation $(a;q)_m:=(1-a)(1-aq)\cdots(1-aq^{m-1}),m\ge 0$ refers to the $q$-pochhammer symbol with $(a;q)=0$.
\end{prop}

\begin{proof}
The decomposition \eqref{eq: decomposition of automorphism} is valid because there are no non-trivial maps between the summands corresponding to each $i$. Also, notice that $\F_p[t]/(f_i)$ is the finite field with $p^{\deg f_i}$ elements. Therefore, the explicit expression \eqref{eq: explicit expression of automorphism} is given in \cite[Chapter 2.1]{macdonald1998symmetric}.
\end{proof}





\section{Proof of the main results}\label{sec: Proof of the main results}

Our starting point is the following observation, formulated as a theorem below. 

\begin{thm}\label{thm: rational canonical and smith normal form}
Let $A_n\in\Mat_n(\F_p)$, and $f\in S$. Then the type of $\Cok(tI_n-A_n)$ at $f$ is $\l^{(f)}=\l^{(f)}(A_n)$.
\end{thm}

\begin{proof}
Notice that the conjugation action of $\GL_n(\F_p)$ does not change the cokernel. Thus, there is no loss of generality in assuming that $A_n$ is already the rational canonical form, i.e.,
$$A_n=\begin{pmatrix}
R_1 & 0 & 0 & \cdots & 0 \\
0 & R_2 & 0 & \cdots & 0 \\
\vdots & \vdots & \vdots & \ddots & \vdots \\
0 & 0 & 0 & \cdots & R_k \\
\end{pmatrix}.$$
Here, each block matrix $R_i$ corresponds to different polynomials in $S$. Denote $\l^{(f)}=(\l_1,\ldots,\l_j)$. We furthermore assume that $R_1$ corresponds to $f$, so that the type of $\Cok(tI
_n-A_n)$ at $f$ is the same as $\Cok(tI_{|\l|\deg f}-R_1)$, and
$$R_1=\begin{pmatrix}
C(f^{\lambda_1}) & 0 & 0 & \cdots & 0 \\
0 & C(f^{\lambda_2}) & 0 & \cdots & 0 \\
\vdots & \vdots & \vdots & \ddots & \vdots \\
0 & 0 & 0 & \cdots & C(f^{\lambda_j}) \\
\end{pmatrix}.$$
In this case, we have
$$\Cok(tI_{|\l|\deg f}-R_1)=\bigoplus_{i=1}^j\Cok(tI_{\l_i\deg f}-C(f^{\l_i}))=\bigoplus_{i=1}^j\left(\F_p[t]/(f^{\lambda_i})\right),$$
which ends the proof.
\end{proof}

The following lemma is a necessary preparation for our proof of \Cref{thm: joint CL distribution}, and may have independent interest.

\begin{lemma}\label{lem: Ext equals Hom}
Let $f\in S$, and $K=\F_p[t]/(f)$ be the corresponding finite field. Let $P\in\F_p[t]$ be monic such that $f\mid P$. Denote by $R:=\F_p[t]/(P)$ the finite quotient ring, and let $M$ be a finite $R$-module. Then, we have
$$\#\Ext^1_R(M,K)=\#\Hom_R(M,K).$$
\end{lemma}

\begin{proof}
Applying \Cref{prop: classification of finite module}, we can construct a short exact sequence of finite $R$-modules
$$0\rightarrow R^s\rightarrow R^s\rightarrow M\rightarrow 0$$
for some $s\in\Z_{\ge 0}$. Using that $\Ext_R^1(R^s,K)=0$, we derive the following exact sequence of $R$-modules:
$$0\rightarrow\Hom_R(M,K)\rightarrow\Hom_R(R_j^s,K)\rightarrow\Hom_R(R_j^s,K)\rightarrow\Ext_R^1(M,K)\rightarrow 0.$$
Therefore, we have
$$\#\Ext^1_R(M,K)=\frac{\#\Hom_R(R_j^s,K)}{\#\Hom_R(R_j^s,K)}\#\Hom_R(M,K)=\#\Hom_R(M,K).$$
\end{proof}

\begin{proof}[Proof of \Cref{thm: joint CL distribution}]
Let $G$ be an arbitrary fixed finite $\F_p[t]$-module. By \cite[Theorem 1.12]{cheong2023distribution}\footnote{In the notation of \cite[Theorem 1.12]{cheong2023distribution}, we are taking $k=1$, $P\in\F_p[t]$ be a monic polynomial such that $PG=0$, and $R=\F_p[t]/(P)$ be the quotient ring. In this case, the surjections of $\F_p[t]$-modules can be naturally regarded as surjections of $R$-modules.}, we have
\begin{equation}\label{eq: expectations are one}
\lim_{n\rightarrow\infty}\E[\#\Sur_{\F_p[t]}(\Cok(tI_n-A_n),G)]=1.
\end{equation}
Now, fix a sequence of partitions $\l^{(1)},\ldots,\l^{(r)}\in\Y$, and denote
$$M:=\bigoplus_{i=1}^r\bigoplus_{j\ge 1}(\F_p[t]/(f_i^{\l^{(i)}_j})),$$
which is a finite $\F_p[t]$-module. Moreover, let 
$$P:=\prod_{i=1}^r f_i^{\lambda_1^{(i)}}\in\F_p[t],$$
so that $PM=0$. Denote by $R:=\F_p[t]/(P)$, so that $M$ can be regarded as a finite $R$-module.
Based on the moments given in \eqref{eq: expectations are one}, we can apply \cite[Lemma 6.3]{sawin2022moment} and \Cref{thm: rational canonical and smith normal form} to deduce that 
\begin{align}\label{eq: from moment to Sawin-Wood}
\begin{split}
\lim_{n\rightarrow\infty}\mathbf{P}(\l^{(f_i)}=\l^{(i)},\forall 1\le i\le r)&=\lim_{n\rightarrow\infty}\mathbf{P}(\text{The type of $\Cok(tI_n-A_n)$ at $f_i$ is $\l^{(i)}$},\forall 1\le i\le r )\\
&=\frac{1}{\Aut(M)}\cdot\prod_{i=1}^r\prod_{j\ge 1}\left(1-\frac{\#\Ext_R^1(M,\F_p[t]/(f_i))}{\#\Hom_{\F_p[t]}(M,\F_p[t]/(f_i))}p^{-j\deg f_i}\right).
\end{split}
\end{align}
Here, let us interpret the notation in \cite[Lemma 6.3]{sawin2022moment}: $u$ is zero, $n$ shall be replaced by $r$, and $N$ shall be replaced by $M$.
For all $1\le i\le r$, $K_i$ is the finite field $\F_p[t]/(f_i)$, thus $q_i=|K_i|$ is equal to $p^{\deg f_i}$. By \Cref{prop: decomposition and explicit expression of Aut}, we have 
$$\#\Aut(M)=\prod_{i=1}^r\#\Aut\left(\bigoplus_{j\ge 1}\F_p[t]/(f_i^{\l^{(i)}_j})\right)=\prod_{i=1}^r\#\Aut_{f_i}(\l^{(i)}).$$
Also, by \Cref{lem: Ext equals Hom}, we have
$$\#\Ext_R^1(M,\F_p[t]/(f_i))=\#\Hom_{\F_p[t]}(M,\F_p[t]/(f_i)),\quad\forall 1\le i\le r.$$
Therefore, we have 
$$\text{RHS}\eqref{eq: from moment to Sawin-Wood}=\prod_{i=1}^r\left(\frac{1}{\#\Aut_{f_i}(\l^{(i)})}\prod_{j\ge 1}(1-p^{-j\deg f_i})\right)=\prod_{i=1}^r\mu_{f_i}(\l^{(i)}).$$
This completes the proof.
\end{proof}

\begin{rmk}[On the relation to Cheong-Yu]
After invoking \Cref{thm: rational canonical and smith normal form} to identify the rational
canonical form data with the $\F_p[t]$-module $\Cok(tI_n-A_n)$, the remainder of our proof
essentially follows the route of Cheong-Yu \cite{cheong2023distribution}. The main
simplification in our finite field setting is that the $\Ext$ and $\Hom$ factors appearing in the general Sawin-Wood formalism collapse, so the expressions reduce to the
clean product forms used in \Cref{thm: joint CL distribution}.

Moreover, Cheong-Yu's \cite[Theorem 1.5]{cheong2023distribution} establishes a universality
statement at the level of $\Z_p[t]$-modules, which tracks strictly more information than the
$\F_p[t]$-module formulation considered here. In particular, once one is aware of
\Cref{thm: rational canonical and smith normal form}, their result should implicitly contain
\Cref{thm: joint CL distribution} as a consequence.
\end{rmk}

Finally, we deduce \Cref{cor: joint multiplicity} from \Cref{thm: joint CL distribution}.

\begin{proof}[Proof of \Cref{cor: joint multiplicity}]
By \cite[Example 5.9(ii)]{cohen2006heuristics}, if $\lambda\in\Y$ is random with respect to the probability measure $\mu_{f_i}$ where $1\le i\le r$,  we have 
$$\mathbf{P}(|\l|=j)=p^{-j\deg f_i}\prod_{k\ge j+1}^\infty(1-p^{-k\deg f_i}),\quad\forall j\in\Z_{\ge 0}.$$
Combining \Cref{thm: joint CL distribution}, this completes the proof.
\end{proof}

\bibliographystyle{plain}
\bibliography{references.bib}

\end{document}